%% file: proth.tex
\documentclass[11pt]{article}
\usepackage{amsmath,amsthm}
\usepackage{amssymb,latexsym}
\usepackage{enumerate}
\usepackage{amstext}
\usepackage{url}

\input{preamble.tex}

\begin{document}
\title{Deterministic Primality Proving on Proth Numbers}
\author{Tsz-Wo Sze ({\tt szetszwo@cs.umd.edu})}
\date{\today}
\maketitle

\renewcommand{\thefootnote}{}


\renewcommand{\thefootnote}{\arabic{footnote}}
\setcounter{footnote}{0}

\begin{abstract}
We present an algorithm to decide the primality of Proth numbers,
$N=2^et+1$,
without assuming any unproven hypothesis.
The expected running time and the worst case running time of the algorithm
are $\tilde{O}((t\log t+\log N)\log N)$
and $\tilde{O}((t\log t+\log N)\log^2 N)$ bit operations,
respectively.
\end{abstract}
\section{Introduction}

A \emph{Proth number} is a positive integer of the form
\begin{equation}\label{eqn-proth numbers}
N=2^e\cdot t+1
\qquad\text{for some odd }t\text{ with }2^e>t>0.
\end{equation}
A self-taught farmer,
Fran\c{c}ois Proth,
published Theorem~\ref{thm-Proth} below in 1878;
see \cite{Williams1998} for more details
and a proof of Theorem~\ref{thm-Proth}.

\begin{theorem}[Proth Theorem]\label{thm-Proth}
Let $N$ be a Proth number defined in (\ref{eqn-proth numbers}).
If
\begin{equation}\label{eqn-proth}
a^{(N-1)/2}\equiv -1\pmod{N}
\end{equation}
for some integer $a$,
then $N$ is a prime.
\end{theorem}

As a consequence,
the primality of Proth numbers can be decided
by a simple, fast probabilistic primality test,
called \emph{the Proth test},
which randomly chooses an integer $a\not\equiv0\pmod{N}$
and then computes 
\begin{equation}\label{eqn-b=a^((N-1)/2) (mod N)}
b\equiv a^{(N-1)/2}\pmod{N}.
\end{equation}
We have the following cases.
\renewcommand{\labelenumi}{(\roman{enumi})}
\begin{enumerate}
\item
If $b\equiv-1\pmod{N}$,
then $N$ is a prime by Theorem~\ref{thm-Proth}.
\item
If $b\not\equiv\pm1\pmod{N}$ and $b^2\equiv 1\pmod{N}$,
then $N$ is composite because $\gcd(b\pm1,N)$ are non-trivial factors of $N$.
\item
If $b^2\not\equiv 1\pmod{N}$,
then $N$ is composite by Fermat's little theorem.
\item
If $b\equiv 1\pmod{N}$,
the primality of $N$ remains unknown.
\end{enumerate}
\renewcommand{\labelenumi}{\arabic{enumi}.}
Case (iv) is the reason that the Proth test is probabilistic.
In case (i),
the integer $a$ is a quadratic nonresidue modulo $N$.
The procedure is repeated
until the primality of $N$ is decided.

When $N$ is a prime,
the Proth test has probability $\frac{1}{2}$ being able to return $N$ prime
in one iteration
since the number of quadratic nonresidues modulo $N$
is exactly $\frac{N-1}{2}$.
The expected number of iterations is $O(1)$.
Equation~(\ref{eqn-b=a^((N-1)/2) (mod N)}) can be computed
in $\tilde{O}(\log^2 N)$ bit operations
using fast integer multiplications \cite{Furer2007, Schoenhage1971}.
Consequently,
the expected running time is $\tilde{O}(\log^2 N)$ bit operations
uniformly for all $t<2^e$.
However,
the worst case running time is exponential
since the number of $a$ such that $b\equiv 1\pmod{N}$ is linear to $N$.

For any integer $N$,
when $N-1$ is fully factored or when $N-1$ is partially factored,
deterministic algorithms by Konyagin and Pomerance can be used
to prove primality in polynomial time \cite{Konyagin1996}.
Their algorithms apply to all Proth numbers $N=2^et+1$,
and run in $\tilde{O}(\log ^{3+C}N)$,
where $C>0$ is a constant depending on $e$ and $t$.
The key idea is the use of \emph{smooth numbers}
so that an ``exponentially large'' subgroup of $(\Z/N\Z)^\times$
can be created by a ``polynomial sized'' set of generators.
We will use smooth numbers to calculate $\sqrt{-1}\pmod{N}$
in Algorithm~\ref{alg-fast sqrt(-1)}.

There are interesting special cases of Proth numbers.
When $t=1$,
the Proth number $N=2^e+1$ is a Fermat number.
It is easy to see that $N$ is composite if $e$ is not a power of two.
In addition,
the primality of $N$ can be decided by \emph{the Pepin test} \cite{Pepin1877},
which is virtually the same as the Proth test,
except that $a=3$ is specified in equation~(\ref{eqn-proth}).
The Pepin test is deterministic due to the following theorem.

\begin{theorem}[Pepin Test]\label{thm-Pepin}
Let $N=2^e+1>3$.
Then $N$ is a prime if and only if
\begin{equation}\label{eqn-pepin}
3^{(N-1)/2}\equiv -1\pmod{N}.
\end{equation}
\end{theorem}

We skip the proof of Theorem~\ref{thm-Pepin}
and show a trivial generalization,
Theorem \ref{thm-Generalized Pepin},
below.
As a result,
for Proth numbers with $t$ not divisible by 3,
the primality can be decided in $\tilde{O}(\log^2 N)$ bit operations.

\begin{theorem}[Generalized Pepin Test]\label{thm-Generalized Pepin}
Let $N=2^e t+1>3$ be a Proth number defined in (\ref{eqn-proth numbers})
such that $3\nmid t$.
Then $N$ is a prime if and only if
\begin{equation}\label{eqn-generalized pepin}
3^{(N-1)/2}\equiv -1\pmod{N}.
\end{equation}
\end{theorem}
\begin{proof}
For any Proth number,
we have $e\geq 1$ by the condition $2^e>t>0$.
If $e=1$,
then $t=1$ by the same condition.
However,
the case $N=3$ is excluded by the assumption $N>3$.
Therefore,
we have $e\geq 2$.

When $N$ is a prime,
we have $3\nmid N$.
Then,
$N\equiv 2\pmod{3}$ by the assumption $3\nmid t$.
The Legendre symbol
$\LEGENDRE{3}{N}=\LEGENDRE{N}{3}=\LEGENDRE{2}{3}=-1$
so that 3 is a quadratic nonresidue modulo $N$.

Finally,
the theorem follows from Theorem~\ref{thm-Proth}.
\end{proof}

Similarly,
there is a deterministic, $\tilde{O}(\log^2 N)$ algorithm
to decide the primality of Proth numbers with $t$ not divisible by 5.

\begin{theorem}\label{thm-5}
Let $N=2^e\cdot t+1>5$ be a Proth number defined in (\ref{eqn-proth numbers})
such that $5\nmid t$.
Then $N$ is a prime if and only if
\[
a^{(N-1)/2}\equiv -1\pmod{N},
\]
where
\begin{equation}\label{eqn-5 not divides t, a}
a=\begin{cases}
5,&\text{if }N\equiv0,2,3\pmod{5}; \\
\frac{-5\pm\sqrt{5}}{2}\pmod{N},&\text{if }N\equiv4\pmod{5}.
\end{cases}
\end{equation}
\end{theorem}
\begin{proof}
Suppose $N$ is a prime with $N\equiv2,3\pmod{5}$.
Then,
$5^{(N-1)/2}\equiv -1\pmod{N}$
since $5$ is a quadratic nonresidue modulo $N$.

Suppose $N$ is a prime with $N\equiv4\pmod{5}$.
Then,
$5$ is a quadratic residue modulo $N$.
Let
\[
\zeta=\frac{\alpha+\sqrt{\alpha^2-4}}{2},
\qquad\text{where }\alpha=\frac{-1\pm\sqrt{5}}{2}\in\frac{\Z}{N\Z}.
\]
Both values of $\zeta$ are 5th roots of unity modulo $N$.
However,
5th roots of unity are not in $\Z/N\Z$
since $N\not\equiv1\pmod{5}$.
Therefore,
$\zeta\not\in\Z/N\Z$
and
$\alpha^2-4=\frac{-5\mp\sqrt{5}}{2}$ are quadratic nonresidues modulo $N$.

The theorem follows from Theorem~\ref{thm-Proth}.
\end{proof}

Note that,
when $N\equiv4\pmod{5}$,
it is unnecessary to calculate $\sqrt{5}\pmod{N}$
in equation (\ref{eqn-5 not divides t, a}).
We may check whether
\[
\left(\frac{-5+x}{2}\right)^{(N-1)/2}\equiv-1\pmod{x^2-5,N}
\]
since,
for $N$ prime,
both $\frac{-5\pm\sqrt{5}}{2}$ are quadratic nonresidues modulo $N$.

Another interesting special case of Proth numbers
is \emph{Cullen numbers} \cite{Cullen1905}.
A Cullen number is a number of the form
\[
C_n=2^n\cdot n+1
\qquad\text{for }n\geq 1.
\]
Notice that $n$ is allowed to be even.
Obviously,
a Cullen number is a Proth number since $2^n>n$ for any positive integer $n$.
There are well-known divisibility properties
of Cullen numbers~\cite{keller1995}.
A list is given below without proofs.
Let $p$ be any odd prime.
\begin{itemize}
\item
Let $n_k=(2^k-k)(p-1)-k$ for $k\geq 0$.
Then, $p\mid C_{n_k}$.
\item
Let $d=\ORD_p 2$.
If $p\mid C_n$,
then $p\mid C_{n+pd}$.
\item
If the Legendre symbol $\LEGENDRE{2}{p}=-1$,
then $p\mid C_{(p+1)/2}$.
\item
If the Legendre symbol $\LEGENDRE{2}{p}=1$,
then $p\mid C_{(3p-1)/2}$.
\end{itemize}

In \cite{stw2011sqrt, stw_phd},
a new idea is introduced for deciding the primality of Proth numbers
but the details are missing.
In this paper,
we fill in all the details
and extend the idea to obtain the following theorem.

\begin{theorem}\label{thm-primality}
Let $N$ be a Proth number defined in (\ref{eqn-proth numbers}).
There is an algorithm deciding the primality of $N$.
The expected running time and the worst case running time of the algorithm are
\[
\tilde{O}((t\log t+\log N)\log N)
\]
and
\[
\tilde{O}((t\log t+\log N)\log^2 N)
\]
bit operations,
respectively.
\end{theorem}

When $t=O(\log N)$,
the expected running time and the worst case running time 
are $\tilde{O}(\log^2 N)$ 
and $\tilde{O}(\log^3 N)$ bit operations,
respectively.
Note that Cullen numbers are covered in this case.
To the best of our knowledge,
our algorithm (Algorithm~\ref{alg-randomized_primality}) is the fastest among
other known \emph{primality proving algorithms}
for numbers of this kind.
By primality proving algorithms,
we mean primality algorithms that always return the correct output.
Adleman-Pomerance-Rumely \cite{Adleman1983} runs in sub-exponential time.
The running times of AKS \cite{Agrawal2004} 
and Lenstra-Pomerance's modified AKS algorithm \cite{Lenstra2009}
are $\tilde{O}(\log^{7.5} N)$ and $\tilde{O}(\log^6 N)$, 
respectively.
The algorithms by Konyagin and Pomerance \cite{Konyagin1996} run
in $\tilde{O}(\log^{3+C} N)$
for some $0<C\leq 3/7$,
$C$ a constant depending on $e$ and $t$.
All the algorithms mentioned above have been proved unconditionally.
With extra assumptions such as the Extended Riemann Hypothesis,
we have the following results:
The elliptic curve primality proving algorithm \cite{Atkin1993, Goldwasser1986}
runs in $\tilde{O}(\log^5N)$.
The running time of Miller's algorithm \cite{Miller1975} is
$\tilde{O}(\log^4N)$.
AKS can be improved \cite{Bernstein2003, Mihailescu2007}
to $\tilde{O}(\log^4 N)$.
The Proth test becomes deterministic
and the running time is $\tilde{O}(\log^4 N)$
since it only has to check congruence equation~(\ref{eqn-proth})
with $2\leq a\leq k$,
where $k$ is $O(\log^2 N)$
by the results from Ankeny \cite{Ankeny1952}.

For $t=O(\log N)$,
Algorithm~\ref{alg-randomized_primality} attains
the same order of expected running time,
$\tilde{O}(\log^2N)$,
as the probabilistic Proth test.
Although Algorithm~\ref{alg-randomized_primality} is a randomized algorithm,
it always returns the correct output in one iteration.
In contrast,
the Proth test
may be unable to decide the primality of $N$
in any fixed number of iterations
because of the exponential time worst cases.
When the actual numbers of bit operations are compared,
the Proth test is faster.
In practice,
a hybrid approach,
which first runs a fixed number of iterations of the Proth test
and then Algorithm~\ref{alg-randomized_primality} if the Proth test fails,
remedies the worst case scenario of the Proth test.

The basic algorithm and the computation of $\sqrt{-1}\pmod{N}$
are described in \SECTION{sect-Primality Proving}.
In \SECTION{sect-Square Root},
we show a square root algorithm
which is the main ingredient of the basic algorithm.
Randomization is described in the Appendix
for improving the expected running time.
\section{Deterministic Primality Proving}\label{sect-Primality Proving}

Let $N>3$ be a Proth number defined in (\ref{eqn-proth numbers}).
We present a deterministic primality proving algorithm,
Algorithm~\ref{alg-primality},
for numbers of this kind in this section.
The proof of correctness and the running time analysis
are shown in the next section.

Fast integer arithmetic is used in all running time analyses.
Therefore,
integer multiplications and divisions modulo $N$ can be computed
in $\tilde{O}(\log N)$ \cite{Knuth1997, Brent2010}.
Denote a fixed square root of $x$ modulo $N$ by $\sqrt{x}\pmod{N}$.

\begin{algorithm}[Deterministic Primality Proving]\label{alg-primality}
The input is $N>3$, a Proth number defined in (\ref{eqn-proth numbers}).
This algorithm returns \PRIME\ if $N$ is a prime.
Otherwise, it returns \COMPOSITE.
\end{algorithm}
\begin{enumerate}
\item[I.]
Try finding $a_2=\sqrt{-1}\pmod{N}$ by Algorithm~\ref{alg-slow sqrt(-1)}. \\
If Algorithm~\ref{alg-slow sqrt(-1)} halts due to $N$ composite,
return \COMPOSITE.
\item[II.]
For each ($3\leq j \leq e$) \{ \\
\hspace*{1cm} Try computing $a_j=\sqrt{a_{j-1}}\pmod{N}$
by Algorithm~\ref{alg-sqrt}. \\
\hspace*{1cm} 
If Algorithm~\ref{alg-sqrt} halts due to $N$ composite, \\
\hspace*{1cm} 
return \COMPOSITE. \\
\}
\item[III.]
Return \PRIME.
\end{enumerate}

We discuss the first step below and the second step,
the crucial step of Algorithm~\ref{alg-primality},
in the next section.

Step~I can be computed as follows:
Suppose $N>3$ is a prime.
We have $e>1$ by definition (\ref{eqn-proth numbers}).
Therefore,
$\sqrt{-1}\in\Z/N\Z$,
where $\Z$ denotes the set of integers
and $\Z/N\Z$ is the prime field with $N$ elements.
Let $H$ be the subgroup of $(\Z/N\Z)^\times$ with $2t$ elements,
where $(\Z/N\Z)^\times$ denotes the multiplicative group of $\Z/N\Z$.
There exists some
$i$, $1\leq i\leq 2t+1$,
such that $i\not\in H$, and
so $i^{2t}\not\equiv 1 \pmod{N}$.
If $i^{2^kt}\equiv-1\pmod{N}$ for some $1\leq k\leq e-1$,
then $i^{2^{k-1}t}\equiv\pm\sqrt{-1}\pmod{N}$.

Suppose $N$ is a prime or a composite number.
If $j^{2t}\equiv1\pmod{N}$ for all $1\leq j\leq 2t+1$,
we deduce that $N$ is composite
since there are $2t+1$ elements with order dividing $2t$.
For some $1\leq i\leq 2t+1$,
if $i^{2t}\not\equiv1\pmod{N}$
but $i^{2^kt}\not\equiv-1\pmod{N}$ for all $1\leq k\leq e-1$,
then either 
\begin{enumerate}
\item[(1)]
$i^{2^et}\not\equiv1\pmod{N}$, or
\item[(2)]
$i^{2^kt}\not\equiv\pm1\pmod{N}$ and $i^{2^{k+1}t}\equiv1\pmod{N}$
for some $1\leq k\leq e-1$.
\end{enumerate}
In case (1),
$N$ is composite by Fermat's little theorem.
In case (2),
$\gcd(i^{2^kt}\pm 1, N)$ are non-trivial factors of $N$ and so $N$ is composite.

\begin{algorithm}[Computing $\sqrt{-1}\pmod{N}$]\label{alg-slow sqrt(-1)}
The input is $N=2^et+1$,
$N$~not necessary a Proth number,
for some integer $e>1$ and odd $t$.
If $N$ is a prime,
this algorithm returns $b$ such that $b^2\equiv-1\pmod{N}$.
Otherwise,
this algorithm either returns an integer $b$ such that $b^2\equiv-1\pmod{N}$
or halts due to $N$ composite.
\end{algorithm}
\renewcommand{\labelenumi}{I.\arabic{enumi}.}
\begin{enumerate}
\item
Compute $b_j=j^{2t}\pmod{N}$ for $1\leq j\leq 2t+1$.
\item
If $b_j\equiv1\pmod{N}$ for all $1\leq j\leq 2t+1$, \\
halt due to $N$ composite.
\item
Suppose $b_i\not\equiv1\pmod{N}$ for some $1\leq i\leq 2t+1$. \\
Compute $b_i^{2^k}\pmod{N}$ for $0\leq k\leq e-2$.
\item
If $b_i^{2^k}\not\equiv-1\pmod{N}$ for all $0\leq k\leq e-2$, \\
halt due to $N$ composite.
\item
Suppose $b_i^{2^{k_0}}\equiv-1\pmod{N}$ for some $0\leq k_0\leq e-2$. \\
Return $i^{2^{k_0}t}\pmod{N}$.
\end{enumerate}

Algorithm~\ref{alg-slow sqrt(-1)} runs in
\[
\tilde{O}((t\log t+\log N)\log N)
\]
bit operations
since steps I.1, I.2 take $\tilde{O}(t\log t\log N)$ bit operations
and steps I.3, I.4, I.5 take $\tilde{O}(e\log N)=\tilde{O}(\log^2 N)$.

For $N$ a Proth number and $t$ large,
the running time can be improved to
\[
\tilde{O}(\log t\log^3 N).
\]
The required $2t+1$ integers in Step~I.1
in Algorithm~\ref{alg-slow sqrt(-1)} can be generated
by a much smaller set.
Let $B$ be a positive integer.
An integer is \emph{$B$-smooth}
if no prime factor of it exceeds $B$.
Let 
\begin{equation}\label{eqn-set of B-smooth integers}
S=
\left\{n\leq N-1\;:\;
n=p_1^{e_1}\cdots p_{\pi(B)}^{e_{\pi(B)}}\text{ with }e_i\geq 0
\right\}
\end{equation}
be the set of positive $B$-smooth integers less than $N$,
where $p_k$ denotes the $k$th prime
and $\pi(\,\cdot\,)$ denotes the prime counting function.

Suppose $N$ is a prime.
Recall that $H$ is the subgroup of $(\Z/N\Z)^\times$ with $2t$ elements.
If $|S|>|H|$,
there exists $i$ such that $p_i\not\in H$
since $p_1,\ldots,p_{\pi(B)}$ generate $S$.
Let $B=\log^A(N-1)$ for some $A>1$.
Then,
$|S|\geq (N-1)^{1-1/A}$ by \cite[Theorem 2.1]{Konyagin1996};
see also \cite[equation (1.14)]{Granville2008}.
Put $A=\frac{2e}{e-1}$.
We have
\[
|S|
\geq (N-1)^{1-\frac{e-1}{2e}}
=(2^et)^{\frac{e+1}{2e}}
=2(2^{\frac{e-1}{2}}t^{\frac{e+1}{2e}})
> 2t.
\]
The last inequality is due to $2^e>t$ in definition (\ref{eqn-proth numbers}).
Algorithm~\ref{alg-slow sqrt(-1)} can be rewritten as below.

\begin{algorithm}[Fast $\sqrt{-1}\pmod{N}$ for $t$ large]
\label{alg-fast sqrt(-1)}
The input is $N>3$, a Proth number defined in (\ref{eqn-proth numbers}).
The output is the same as the one in Algorithm~\ref{alg-slow sqrt(-1)}.
\end{algorithm}
\renewcommand{\labelenumi}{I.\alph{enumi}.}
\begin{enumerate}
\item
Let $B=\lfloor\log^{2e/(e-1)}(N-1)\rfloor$. \\
Find $p_1,\ldots,p_{\pi(B)}$, the primes less than or equal to $B$.
\item
Compute $b_j=p_j^{2t}\pmod{N}$ for $1\leq j\leq \pi(B)$.
\item
If $b_j\equiv1\pmod{N}$ for all $1\leq j\leq \pi(B)$, \\
halt due to $N$ composite.
\item
Suppose $b_i\not\equiv1\pmod{N}$ for some $1\leq i\leq \pi(B)$. \\
Compute $b_i^{2^k}\pmod{N}$ for $0\leq k\leq e-2$.
\item
If $b_i^{2^k}\not\equiv-1\pmod{N}$ for all $0\leq k\leq e-2$, \\
halt due to $N$ composite.
\item
Suppose $b_i^{2^{k_0}}\equiv-1\pmod{N}$ for some $0\leq k_0\leq e-2$. \\
Return $p_i^{2^{k_0}t}\pmod{N}$.
\end{enumerate}

Step~I.a takes $\tilde{O}(\log^{2+2/(e-1)} N)$
since 
calculating $B$ requires $\tilde{O}(\log N)$
and the primes $p_1,\ldots,p_{\pi(B)}$ can be obtained
in $\tilde{O}(B)$ by a sieve \cite{Atkin2004}.
The remaining steps are essentially
the same as Algorithm~\ref{alg-slow sqrt(-1)}
except that number of iterations is decreased
from $2t+1$ in Step~I.1 to $\pi(B)$ in Step~I.b.
The running time of the steps from I.b to I.f is
$\tilde{O}(\log t\log^{3+2/(e-1)} N)$,
which is also the overall running time of Algorithm~\ref{alg-fast sqrt(-1)}.
Finally,
$O(\log^{2/(e-1)} N)=O(1)$
because $N$ is a Proth number.
\section{Taking Square Roots}\label{sect-Square Root}

Suppose $N$ is a prime,
but not necessary a Proth prime,
with $N\equiv1\pmod{4}$ for the following.
Given $\sqrt{-1}\pmod{N}$,
a square root of a fixed value,
we show how to calculate the square roots of an arbitrary value $\beta$
for $1<\beta<N-1$ and $\beta$ a quadratic residue modulo $N$.
Such ideas of taking square roots modulo $N$ are presented
in \cite{stw2011sqrt, stw_phd}.
We include the material below for completeness.
Additionally,
we show how to detect the case that $N$ is a composite number.

Suppose
\[
\beta\equiv\alpha^2\pmod N
\qquad\text{for some integer } \alpha.
\]
A group $\G$ is constructed
such that $\G\simeq(\Z/N\Z)^\times$.
The isomorphism $\psi_\alpha:\G\longrightarrow(\Z/N\Z)^\times$
depends on $\alpha$ as a parameter.
Then,
we find an order 4 element in $\G$.
Such element must be mapped to $\pm\sqrt{-1}\pmod{N}$ through $\psi_\alpha$.
Consequently,
the parameter $\alpha$ of $\psi_\alpha$ can be calculated.

It can be shown that the group $\G$ is isomorphic to
a ``singular elliptic curve''.
Coincidentally,
elliptic curve arithmetic is also used
in Schoof's square root algorithm over prime fields \cite{Schoof1985}.
However,
the ideas of these two algorithms are quite different.
In Schoof's algorithm,
an elliptic curve $E$ is constructed
over a selected extension of $\Z/N\Z$ for $N$ prime
in such a way that $E$ has complex multiplication.
Then,
compute the Frobenius endomorphism
and,
subsequently,
obtain the square roots modulo $N$.

We begin describing our algorithm by
defining two sets,
\begin{eqnarray}
\G'&\DEF{=}& \SETR{\ELEMENT{a}}{a\not\equiv\pm\alpha\pmod{N}},\text{ and}\\
\G &\DEF{=}& \G'\cup\SET{\ELEMENT{\infty}}.
\end{eqnarray}
The elements in $\G$ are denoted by $\ELEMENT{\,\cdot\,}$
for avoiding confusion with the elements in ${\cal Z}$,
where
\begin{eqnarray}\label{eqn-Z}
{\cal Z} &\DEF{=}& \Z\cup\SET{\infty}.
\end{eqnarray}
Two elements $\ELEMENT{a_1}, \ELEMENT{a_2}\in\G'$ are equal
if and only if $a_1\equiv a_2\pmod{N}$.
Therefore,
there are $N-2$ and $N-1$ elements in $\G'$ and $\G$,
respectively.

Further,
define an operator $*$ as follows:
For any $\ELEMENT{a}\in\G$
and any $\ELEMENT{a_1}, \ELEMENT{a_2}\in\G'$ with $a_1+a_2\not\equiv0\pmod{N}$,
\begin{eqnarray}
\label{eqn-a*infty}
\ELEMENT{a}*\ELEMENT{\infty}&=&\ELEMENT{\infty}*\ELEMENT{a}=\ELEMENT{a}, \\
\label{eqn-a*-a}
\ELEMENT{a_1}*\ELEMENT{-a_1}&=&\ELEMENT{\infty}, \\
\label{eqn-a_1*a_2}
\ELEMENT{a_1}*\ELEMENT{a_2}&=&\ELEMENT{(a_1a_2+\alpha^2)(a_1+a_2)^{-1}},
\end{eqnarray}
where $x^{-1}$ denotes the multiplicative inverse of $x\pmod{N}$
for integer $x$ with $\gcd(x,N)=1$.
We have the following lemma.

\begin{lemma}\label{lem-G=F_N'}
When $N$ is prime,
$(\G,*)$ is a well-defined group,
which is isomorphic to $(\Z/N\Z)^\times$.
\end{lemma}
\begin{proof}
Define a bijective mapping
\begin{equation}\label{eqn-psi}
\psi:\G\longrightarrow\left(\frac{\Z}{N\Z}\right)^\times,
\qquad \ELEMENT{\infty}\longmapsto 1,
\qquad \ELEMENT{a}\longmapsto (a+\alpha)(a-\alpha)^{-1}
\end{equation}\label{eqn-psi^-1}
with inverse mapping
\begin{equation}
\psi^{-1}:\left(\frac{\Z}{N\Z}\right)^\times\longrightarrow\G,
\qquad 1\longmapsto \ELEMENT{\infty},
\qquad b\longmapsto \ELEMENT{\alpha(b+1)(b-1)^{-1}}.
\end{equation}
A straightforward calculation shows that $\psi$ is a homomorphism.
\end{proof}

Note that $\G$ is also isomorphic to the group of non-singular points
$(x,y)\in(\Z/N\Z)^2$
on the ``singular elliptic curve''
\[
y^2=x^2(x+\alpha^2).
\]
For more details,
see \cite{lcw2008}.

In the rest of the paper,
we drop the symbol $*$
and denote the group operation of $\G$ by multiplication.
Algorithm~\ref{alg-multiplication} below shows
how to perform the group operation.
The input integer $N$ may be a prime or a composite number
since the algorithm will be used for deciding the primality of $N$.

\begin{algorithm}[Group Operation]\label{alg-multiplication}
The inputs are $N,\beta\in\Z$ and $a_1,a_2\in{\cal Z}$ such that
$0<\beta<N$
and either $a_i=\infty$ or $a_i^2\not\equiv\beta\pmod{N}$ for $i=1,2$.
If $N$ is a prime, 
the input $\beta$ is guaranteed to be a quadratic residue modulo $N$
and this algorithm returns $a\in{\cal Z}$
such that $\ELEMENT{a}=\ELEMENT{a_1}\ELEMENT{a_2}\in\G$
for $\alpha^2\equiv\beta\pmod{N}$.
Otherwise,
this algorithm either returns some $a'\in{\cal Z}$
or halts due to $N$ composite.
\end{algorithm}
\renewcommand{\labelenumi}{\arabic{enumi}.}
\begin{enumerate}
\item If $a_1=\infty$, return $a_2$.
\item If $a_2=\infty$, return $a_1$.
\item If $a_1+a_2\equiv0\pmod{N}$, return $\infty$.
\item If $\gcd(a_1+a_2, N)\neq1$,
halt due to $N$ composite.
\item Compute $a\equiv(a_1a_2+\beta)(a_1+a_2)^{-1}\pmod{N}$.
\item
If $a^2\equiv\beta\pmod{N}$,
halt due to $N$ composite. \\
Otherwise,
return $a$.
\end{enumerate}

Algorithm~\ref{alg-multiplication} basically follows
definitions~(\ref{eqn-a*infty}), (\ref{eqn-a*-a}) and (\ref{eqn-a_1*a_2}).
It also handles the case if $N$ is a composite number.
In such case,
$\G$ is no longer a well-defined group.
If the algorithm halts in Step~4,
a non-trivial factor of $N$ is discovered and so $N$ is a composite number.
If it halts in Step~6,
we have $(a_1a_2+\beta)^2\equiv\beta(a_1+a_2)^2\pmod{N}$,
which implies $(a_1^2-\beta)(a_2^2-\beta)\equiv0\pmod{N}$.
Since $a_i^2\not\equiv\beta\pmod{N}$ for $i=1,2$ by the assumption,
$a_1^2-\beta$ and $a_2^2-\beta$ are zero divisors,
which means that $N$ is composite.
Note that the value of $\alpha$ is not required
in Algorithm~\ref{alg-multiplication}.

Equipped with $\G$,
we are ready to describe the square root algorithm,
which is the main ingredient of the Step~II in Algorithm~\ref{alg-primality}.
We continue with the following lemma.

\begin{lemma}\label{lem-alpha=a b}
Let $N\equiv1\pmod{4}$ be a prime.
If $\ELEMENT{a}\in\G$ is an order 4 element
and $b\equiv\pm\sqrt{-1}\pmod{N}$,
an order 4 element in $(\Z/N\Z)^\times$,
then
\[
\alpha=\pm a b\pmod{N}.
\]
\end{lemma}
\begin{proof}
By Lemma~\ref{lem-G=F_N'},
we have $\ELEMENT{a}=\psi^{-1}(\pm b)$.
Then,
\[
\alpha
\equiv\frac{a(\pm b-1)}{(\pm b+1)}
\equiv\pm a b
\pmod{N}.
\]
The lemma follows.
\end{proof}

We have Algorithm~\ref{alg-sqrt} below for taking square roots modulo $N$.
The notation $\ELEMENT{x}^y$ means
using Algorithm~\ref{alg-multiplication} and the successive squaring method
to compute $\ELEMENT{x}$ to the power $y$.
Algorithm~\ref{alg-sqrt} halts due to $N$ composite 
as soon as Algorithm~\ref{alg-multiplication} does,
if it is the case.

\begin{algorithm}[Taking Square Root Modulo $N$]\label{alg-sqrt}
The inputs are integers $N$, $\beta$ and $b$
such that $1<\beta<N-1$ and $b^2\equiv-1\pmod{N}$,
where $N=2^et+1$ for some integer $e>1$ and odd $t$.
If $N$ is a prime, 
the input $\beta$ is guaranteed to be a quadratic residue modulo $N$
and this algorithm returns $\alpha$ such that $\alpha^2\equiv\beta\pmod{N}$.
If $N$ is a composite number,
this algorithm either returns an integer $\alpha$
such that $\alpha^2\equiv\beta\pmod{N}$
or halts due to $N$ composite.
\end{algorithm}
\renewcommand{\labelenumi}{II.\arabic{enumi}.}
\begin{enumerate}
\item
Check easy cases:
If $j^2\equiv\beta\pmod{N}$ for $1\leq j\leq 2t+1$, return $j$.
\item
Find $\ELEMENT{a}$ such that 
$\ELEMENT{a}^2\neq\ELEMENT{\infty}$ and $\ELEMENT{a}^4=\ELEMENT{\infty}$
as below:
\begin{enumerate}
\item
Compute $\ELEMENT{c_j}=\ELEMENT{j}^{2t}$ for $1\leq j\leq 2t+1$.
\item
If $\ELEMENT{c_j}=\ELEMENT{\infty}$ for all $1\leq j\leq 2t+1$,
halt due to $N$ composite.
\item
Suppose $\ELEMENT{c_i}\neq\ELEMENT{\infty}$ for some $1\leq i\leq 2t+1$. \\
Compute $\ELEMENT{c_i}^{2^k}$ for $0\leq k\leq e-2$.
\item
If $\ELEMENT{c_i}^{2^k}\neq\ELEMENT{0}$ for all $0\leq k\leq e-2$,
halt due to $N$ composite.
\item
Suppose $\ELEMENT{c_i}^{2^{k_0}}=\ELEMENT{0}$ for some $0\leq k_0\leq e-2$. \\
Compute $\ELEMENT{a}=\ELEMENT{i}^{2^{k_0}t}$.
\end{enumerate}
\item
Compute $\alpha$:
\begin{enumerate}
\item
Compute $\alpha\equiv ab\pmod{N}$.
\item
If $\alpha^2\not\equiv\beta\pmod{N}$,
halt due to $N$ composite. \\
Otherwise,
return $\alpha$.
\end{enumerate}
\end{enumerate}

\begin{proposition}
Algorithm~\ref{alg-sqrt} is correct.
\end{proposition}
\begin{proof}
Step~II.1 checks some easy cases.
If $j$ is a square root of $\beta$ modulo $N$ for some $1\leq j\leq 2t+1$,
we are done.

Step~II.2 is similar to Algorithm~\ref{alg-slow sqrt(-1)}.
For $N$ prime,
both of them compute an order 4 element:
Algorithm~\ref{alg-slow sqrt(-1)} calculates 
an element of $(\Z/N\Z)^\times$ of order 4,
while an element of $\G$ of order 4 is computed in Step~II.2.
Note that $\G$ is isomorphic to the cyclic group $(\Z/N\Z)^\times$ as shown
in Lemma~\ref{lem-G=F_N'}.
The identity element and the order 2 element of $\G$ 
are $\ELEMENT{\infty}$ and $\ELEMENT{0}$,
respectively.

Suppose $N$ is a prime.
Let
\[
H_\alpha\DEF{=}
\SETR{\ELEMENT{g}\in\G}{\ELEMENT{g}^{2t}=\ELEMENT{\infty}},
\]
be the $2t$-torsion subgroup of $\G$.
The size of $H_\alpha$ is $2t$ since $\G$ is cyclic.
If $\ELEMENT{c_j}=\ELEMENT{j}^{2t}=\ELEMENT{\infty}$ for all $1\leq j\leq 2t+1$,
there are $2t+1$ elements with order dividing $2t$ in $\G$,
which leads to a contradiction.
Suppose $\ELEMENT{c_i}\neq\ELEMENT{\infty}$ for some $1\leq i\leq 2t+1$.
In Step~II.2.(d),
if $\ELEMENT{c_i}^{2^k}\neq\ELEMENT{0}$ for all $0\leq k\leq e-2$,
then $\ELEMENT{i}^{2^{e-1}t}\neq\ELEMENT{0}$,
which implies $\ELEMENT{i}^{|\G|}\neq\ELEMENT{\infty}$,
a contradiction.
Thus,
if the algorithm halts at Step~II.2.(b) or Step~II.2.(d),
then $N$ is composite.
The order of the element $\ELEMENT{a}\in\G$ obtained in Step~II.2.(e)
is exactly $4$
since $\ELEMENT{a}^2=\ELEMENT{0}\neq\ELEMENT{\infty}$
and $\ELEMENT{a}^4=\ELEMENT{\infty}$.
By Lemma~\ref{lem-alpha=a b},
$\alpha=a b\pmod{N}$ is square root of $\beta$.
If $\alpha^2\not\equiv\beta\pmod{N}$,
we have $N$ a composite number.

The proposition follows.
\end{proof}

\begin{proposition}
Algorithm~\ref{alg-sqrt} runs in
\[
\tilde{O}((t\log t+\log N)\log N)
\]
bit operations.
\end{proposition}
\begin{proof}
All the powers $\ELEMENT{x}^y$
are computed by Algorithm~\ref{alg-multiplication} 
and the successive squaring method,
which take $\tilde{O}(\log y\log N)$ bit operations.

Step~II.1 takes $\tilde{O}(t\log N)$ bit operations.
In Step~II.2,
the running time of parts (a) and (b) together is
$\tilde{O}(t\log t\log N)$,
parts (c) and (d) together take
$\tilde{O}(e\log N)=\tilde{O}(\log^2 N)$ bit operations,
and part (e) takes $\tilde{O}(\log^2 N)$ bit operations.
Therefore,
Step~II.2 takes $\tilde{O}((t\log t+\log N)\log N)$ bit operations in total.
Step~II.3 only takes $\tilde{O}(\log N)$ bit operations.

The proposition follows.
\end{proof}

Finally,
we show the following propositions.

\begin{proposition}
Algorithm~\ref{alg-primality} is correct.
\end{proposition}
\begin{proof}
If either Step~I or Step~II in Algorithm~\ref{alg-primality} returns \COMPOSITE,
the input $N$ must be \COMPOSITE\ 
by Algorithm~\ref{alg-slow sqrt(-1)} and Algorithm~\ref{alg-sqrt}.
Note that, for $N$ prime,
the integer $a_{j-1}$ in Step~II is a quadratic residue modulo $N$
for all $3\leq j\leq e$.
Otherwise,
Step~III returns \PRIME.
In this case,
$N$ is indeed a prime by Theorem~\ref{thm-Proth} with $a=a_e$,
where $a_e$ is computed in Step~II.
The proposition follows.
\end{proof}

\begin{proposition}
Algorithm~\ref{alg-primality} runs in 
\[
\tilde{O}((t\log t+\log N)\log^2 N)
\]
bit operations.
\end{proposition}
\begin{proof}
Step~I takes $\tilde{O}((t\log t+\log N)\log N)$
bit operations
using Algorithm~\ref{alg-slow sqrt(-1)}.
With Algorithm~\ref{alg-sqrt},
Step~II takes 
$\tilde{O}(e(t\log t+\log N)\log N)=\tilde{O}((t\log t+\log N)\log^2 N)$
bit operations.
Step~III can be done in $O(1)$.
The proposition follows.
\end{proof}

It is not hard to see that
the expected running time of Algorithm~\ref{alg-primality} can be improved to
\[
\tilde{O}((t\log t+\log N)\log N)
\]
by randomizing Step I.
The details are presented in the Appendix.
\begin{proof}[Proof of Theorem~\ref{thm-primality}]
It follows from Propositions~\ref{prop-randomized_primality correctness}
and \ref{prop-randomized_primality running time}.
\end{proof}
\section{Appendix: Randomization}

In this section,
randomization is introduced for improving the expected running time.
Unlike some probabilistic algorithms,
e.g.~the Proth test,
which may not able to decide the primality of $N$,
our randomized algorithm always returns the correct output.
Algorithm~\ref{alg-primality} first tries computing
$a_2\equiv\sqrt{-1}\pmod{N}$,
and then it repeatedly takes square roots to obtain
$a_3$, $a_4$, $\ldots$, $a_e$
such that $a_j\equiv\sqrt{a_{j-1}}\pmod{N}$ for $3\leq j\leq e$.
If $N$ is a prime,
all the computations succeed
and it ends up with $a_e$,
a quadratic nonresidue modulo $N$.
In  total,
this uses $e-1=O(\log N)$ square root computations,
which dominate the running time of the entire algorithm.
The expected running time of Algorithm~\ref{alg-primality} can be improved
by repeatedly taking square roots on a randomly chosen integer,
instead of  the fixed integer $-1$.
We first randomly choose an integer $a$.
Then,
we compute $\sqrt{a}\pmod{N}$, $\sqrt{\sqrt{a}}\pmod{N}$
and so on.
If $N$ is a prime,
this process ends up with a quadratic nonresidue modulo $N$.

For prime $N=2^et+1$,
the multiplicative group $(\Z/N\Z)^\times$ being cyclic tells that 
most of elements in $(\Z/N\Z)^\times$ have order with large 2-part.
Only a few square root computations are required
in order to obtain a quadratic nonresidue from these elements.
In fact,
half of the total number of elements in $(\Z/N\Z)^\times$
are quadratic nonresidues modulo $N$.
The order of a quadratic nonresidue is divisible by $2^e$.
In general,
for $1\leq k\leq e$,
there are exactly $2^{k-1}t$ elements having order divisible by $2^k$
but not $2^{k+1}$.
Only $e-k$ square root computations are required
for obtaining a quadratic nonresidue
from these $2^{k-1}t$ elements.

The randomized algorithm is presented below.

\begin{algorithm}[Randomized Algorithm]\label{alg-randomized_primality}
The input is $N>3$, a Proth number defined in (\ref{eqn-proth numbers}).
This algorithm returns \PRIME\ if $N$ is a prime.
Otherwise, it returns \COMPOSITE.
\end{algorithm}
\renewcommand{\labelenumi}{(\roman{enumi})}
\begin{enumerate}
\item
Find $b_k$ such that the order $b_k\pmod{N}$ is $2^k$
for some $k\geq 2$ as below:
\begin{enumerate}
\item
Randomly choose an integer $1<a<N-1$ \\
until $a^{2t}\not\equiv 1\pmod{N}$. \\
If there are $2t-1$ distinct integers $1<a<N-1$ \\
such that $a^{2t}\equiv 1\pmod{N}$, return \COMPOSITE.
\item
Compute $a_0\equiv a^t\pmod{N}$. \\
If $a_0^{2^e}\not\equiv 1\pmod{N}$, return \COMPOSITE.
\item
Find the least $k\geq 2$ such that $a_0^{2^k}\equiv 1\pmod{N}$. \\
If $a_0^{2^{k-1}}\not\equiv -1\pmod{N}$, return \COMPOSITE.
\item
Set $b_2\equiv a_0^{2^{k-2}}\pmod{N}$. \\
Set $b_k=a_0$.
\end{enumerate}
\item
For each ($k+1\leq j \leq e$) \{ \\
\hspace*{1cm} Try computing $b_j=\sqrt{b_{j-1}}\pmod{N}$
by Algorithm~\ref{alg-sqrt}. \\
\hspace*{1cm} 
If Algorithm~\ref{alg-sqrt} halts due to $N$ composite, \\
\hspace*{1cm} 
return \COMPOSITE. \\
\}
\item
Return \PRIME.
\end{enumerate}

\begin{proposition}\label{prop-randomized_primality correctness}
Algorithm~\ref{alg-randomized_primality} is correct.
\end{proposition}
\begin{proof}
In Step~(i)(a),
we randomly choose an integer $a$ from the open interval $(1,N-1)$
without replacement
until $a^{2t}\not\equiv 1\pmod{N}$.
If there are $2t-1$ distinct integers $a$ in $(1,N-1)$ 
such that $a^{2t}\equiv 1\pmod{N}$,
then these $2t-1$ distinct integers together with $1$ and $N-1$
are totally $2t+1$ distinct integers with order modulo $N$ dividing $2t$.
Therefore,
$N$ is composite.
In Step~(i)(b),
if $a_0^{2^e}\equiv a^{N-1}\not\equiv 1\pmod{N}$,
then $N$ is composite by Fermat's little theorem.
In Step~(i)(c),
the least positive integer $k$ with $a_0^{2^k}\equiv 1\pmod{N}$ exists
since $a_0^{2^e}\equiv 1\pmod{N}$.
We also have $k\geq 2$ because $a^{2t}\not\equiv 1\pmod{N}$ by Step~(i)(a).
If $a_0^{2^{k-1}}\not\equiv -1\pmod{N}$,
then $\gcd(a_0^{2^{k-1}}-1,N)$ is a non-trivial factor of $N$,
and so $N$ is composite.
If Step~(i)(a), (i)(b) and (i)(c) do not return \COMPOSITE,
we end up in Step~(i)(d) that
$b_2\equiv a_0^{2^{k-2}}\equiv\pm\sqrt{-1}\pmod{N}$
and $b_k=a_0$ with $b_k^{2^{k-1}}\equiv-1\pmod{N}$.
Note that $b_2$ is required
as an input of Algorithm~\ref{alg-sqrt} used in Step~(ii).
We will show that the value of $k$ is large with high probability
later in the section.

Step~(ii) and (iii) are similar to the Step~II and III
in Algorithm~\ref{alg-primality}
except that Step~(ii) begins taking square roots with $b_k$.
If Algorithm~\ref{alg-sqrt} does not halt due to $N$ composite,
the invariant $b_j^{2^{j-1}}\equiv-1\pmod{N}$ is maintained in the loop
for $k\leq j\leq e$.
If $b_e$ is obtained,
$N$ is a prime by Theorem~\ref{thm-Proth} with $a=b_e$.

The proposition follows.
\end{proof}

\begin{proposition}\label{prop-randomized_primality running time}
The expected running time and the worst case running time 
of Algorithm~\ref{alg-randomized_primality} 
are
\[
\tilde{O}((t\log t+\log N)\log N)
\]
and
\[
\tilde{O}((t\log t+\log N)\log^2 N)
\]
bit operations,
respectively.
\end{proposition}
\begin{proof}
Step~(i)(a) requires $\tilde{O}(t\log t\log N)$ bit operations.
Steps (i)(b), (i)(c) and (i)(d) together take $\tilde{O}(\log^2 N)$
bit operations.
The running time of Step~(ii) is
$\tilde{O}(m(t\log t+\log N)\log N)$ bit operations,
where $m$ is an upper bound of the number of iterations in the loop.
Step~(iii) can be done in $\tilde{O}(\log N)$ bit operations.
The entire algorithm is dominated by Step~(ii).
The total running time is
$\tilde{O}(m(t\log t+\log N)\log N)$ bit operations.

The value of $m$ depends on the integer $a$ chosen in Step~(i)(a).
It is easy to see that the worst case is $m=O(\log N)$.
We will show that the expected value of $m$ is less than 1
in Lemma~\ref{lem-E(m)=O(1)}.
The proposition follows.
\end{proof}

Let $v_2(x)$ be the 2-adic valuation function.
For positive integer $x=2^rs$ with $s$ odd, 
we have $v_2(x)=r$,
which is the exponent of the 2-part of $x$.
Let $\ORD_p a$ be the order of $a\pmod{p}$
for prime $p$ and $a\not\equiv 0\pmod{p}$.
We show in Lemma~\ref{lem-E(v_2(ord_p a))} below that
the expected value of $v_2(\ORD_p a)$
for a random integer $a$ is bounded below by $v_2(p-1)-1$.

\begin{lemma}\label{lem-E(v_2(ord_p a))}
Let $p=2^{e'}t'+1$ be an odd prime for some odd $t'$ and $e'\geq 1$.
Let $a$ be an integer randomly chosen from the open interval $(1,p-1)$
such that $a^{2d}\not\equiv 1\pmod{p}$
for some positive divisor $d$ of $t'$.
Then the expected value 
\begin{equation*}
E(v_2(\ORD_p a))>e'-1.
\end{equation*}
\end{lemma}
\begin{proof}
By counting the number of integers $a\in(1,p-1)$ such that $v_2(\ORD_p a)=i$
for $i=0,1,\cdots,e'$,
we have
\begin{eqnarray*}
\sum_{\substack{1<a<p-1 \\ a^{2d}\not\equiv 1\pmod{p}}}v_2(\ORD_p a)
&=&  0\cdot (t'-d) + 1\cdot (t'-d) + \sum_{i=2}^{e'}i\cdot 2^{i-1}t' \\
&=& (e'-1)(p-1)+t'-d.
\end{eqnarray*}
Then, the expected value is
\[
E(v_2(\ORD_p a))
= e'-1+\frac{2d(e'-1)+t'-d}{p-1-2d} 
> e'-1.
\]
The lemma follows.
\end{proof}

\begin{lemma}\label{lem-E(m)=O(1)}
\[
E(m)<1.
\]
\end{lemma}
\begin{proof}
Suppose $N$ is a prime.
Recall that $a$ is a randomly chosen integer in Step~(i)(a)
such that $a^{2t}\not\equiv1\pmod{N}$
and $k=v_2(\ORD_N a)$.
By Lemma~\ref{lem-E(v_2(ord_p a))} with $d=t$,
we have $E(k) = E(v_2(\ORD_N a)) > e-1$.
Therefore,
\[
E(m)=E(e-k)<1.
\]

Suppose $N$ is composite.
Let $p$ be a prime divisor of $N$
such that $v_2(p-1)$ is the minimum among all the prime divisors of $N$.
Write $p=2^{e'}t'+1$.
Clearly,
we have $e'\leq e$.
If the algorithm does not discover $N$ composite in Step~(i),
the maximum number of iterations is bounded above by $e'$,
i.e.\ $m\leq e'$.
Let $a$ be the integer chosen in Step~(i)(a).
If $p$ divides $a$,
then Step~(i)(b) will return \COMPOSITE\ 
since $a_0^{2^e}\equiv a^{2^et}\not\equiv 1\pmod{N}$ .
Suppose $p$ does not divide $a$.
By Lemma~\ref{lem-E(v_2(ord_p a))} with $d=\gcd(t,t')$,
we have $E(v_2(\ORD_p a)) > e'-1$.
Finally,
\[
E(m)\leq E(e'-v_2(\ORD_p a))<1.
\]

The lemma follows.
\end{proof}
\bibliographystyle{plain}
\bibstyle{plain}
\bibliography{proth}
\end{document}

%% file: preamble.tex
\newcounter{defcounter}
\newtheorem{theorem}{Theorem}[section]

\newtheorem{lemma}[theorem]{Lemma}
\newtheorem{proposition}[theorem]{Proposition}

\newtheorem{algorithm}[theorem]{Algorithm}

\numberwithin{equation}{section}

\newcommand{\IGNORE}[1]{}
\newcommand{\SECTION}[1]{\S\ref{#1}}
\newcommand{\DEF}[1]{\stackrel{\sf def}{#1}}

\def\Z{{\mathbb Z}}

\def\G{G_\alpha}

\DeclareMathOperator{\ORD}{ord}

\def\PRIME{{\sf PRIME}}
\def\COMPOSITE{{\sf COMPOSITE}}

\newcommand{\ELEMENT}[1]{\left[ #1 \right]}

\newcommand{\LEGENDRE}[2]{\genfrac{(}{)}{}{}{#1}{#2}}

\newcommand{\SET}[1]{\left\{ #1 \right\}}

\newcommand{\SETR}[2]{\left\{ #1 \;:\; #2 \right\}}

%% file: proth.bbl
\begin{thebibliography}{10}

\bibitem{Adleman1983}
Leonard~M. Adleman, Carl Pomerance, and Robert~S. Rumely.
\newblock On distinguishing prime numbers from composite numbers.
\newblock {\em Ann. of Math.}, 117(1):173--206, 1983.

\bibitem{Agrawal2004}
Manindra Agrawal, Neeraj Kayal, and Nitin Saxena.
\newblock {PRIMES} is in {P}.
\newblock {\em Ann. of Math.}, 160(2):781--793, 2004.

\bibitem{Ankeny1952}
Nesmith~C. Ankeny.
\newblock The least quadratic non residue.
\newblock {\em Ann. of Math.}, 55(1):65--72, 1952.

\bibitem{Atkin2004}
A.~O.~L. Atkin and Daniel~J. Bernstein.
\newblock Prime sieves using binary quadratic forms.
\newblock {\em Math. Comp.}, 73(246):1023--1030, 2004.

\bibitem{Atkin1993}
A.~O.~L. Atkin and Fran\c{c}ois Morain.
\newblock Elliptic curves and primality proving.
\newblock {\em Math. Comp.}, 61(203):29--68, 1993.

\bibitem{Bernstein2003}
Daniel~J. Bernstein.
\newblock Proving primality in essentially quartic random time.
\newblock {\em Math. Comp.}, 76(257):389--403, 2007.

\bibitem{Brent2010}
Richard~P. Brent and Paul Zimmermann.
\newblock {\em Modern Computer Arithmetic}.
\newblock Number~18 in Cambridge Monographs on Computational and Applied
  Mathematics. Cambridge University Press, Cambridge, United Kingdom, 2010.

\bibitem{Cullen1905}
James Cullen.
\newblock Question 15897.
\newblock {\em Educ. Times}, page 534, December 1905.

\bibitem{Furer2007}
Martin F{\"u}rer.
\newblock Faster integer multiplication.
\newblock In {\em Proceedings of the 39th Annual ACM Symposium on Theory of
  Computing}, pages 57--66. ACM, 2007.

\bibitem{Goldwasser1986}
Shafi Goldwasser and Joe Kilian.
\newblock Almost all primes can be quickly certified.
\newblock In {\em Proceedings of the 18th Annual ACM Symposium on Theory of
  Computing}, pages 316--329. ACM, 1986.

\bibitem{Granville2008}
Andrew Granville.
\newblock Smooth numbers: computational number theory and beyond.
\newblock In Joseph~P. Buhler and Peter Stevenhagen, editors, {\em Algorithmic
  Number Theory: Lattices, Number Fields, Curves and Cryptography}, number~44
  in Mathematical Sciences Research Institute Publications, pages 267--323.
  Cambridge University Press, 2008.

\bibitem{Lenstra2009}
Hendrik W.~Lenstra Jr. and Carl Pomerance.
\newblock Primality testing with {Gaussian} periods, 2009.
\newblock Preprint (\url{http://math.dartmouth.edu/~carlp/aks102309.pdf}).

\bibitem{keller1995}
Wilfrid Keller.
\newblock New {Cullen} primes.
\newblock {\em Math. Comp.}, 64(212):1733--1741, 1995.

\bibitem{Knuth1997}
Donald~E. Knuth.
\newblock {\em The Art of Computer Programming, Volume 2: Seminumerical
  Algorithms}.
\newblock Addison-Wesley, Reading, MA, 1997.

\bibitem{Konyagin1996}
S.~Konyagin and C.~Pomerance.
\newblock On primes recognizable in deterministic polynomial time.
\newblock In {\em The Mathematics of Paul Erd{\"o}s}, volume~13 of {\em
  Algorithms Combin.}, pages 176--198. Springer-Verlag, Berlin, 1996.

\bibitem{Mihailescu2007}
Preda {Mih\u ailescu} and Roberto~M. Avanzi.
\newblock Efficient ``quasi''-deterministic primality test improving {AKS},
  2007.
\newblock Preprint
  (\url{http://www.uni-math.gwdg.de/preda/mihailescu-papers/ouraks3.pdf}).

\bibitem{Miller1975}
Gary~L. Miller.
\newblock Riemann's hypothesis and tests for primality.
\newblock In {\em Proceedings of Seventh Annual Symposium on Theory of
  Computing}, pages 234--239. ACM, 1975.

\bibitem{Pepin1877}
P.~P\'{e}pin.
\newblock Sur la formule $2^{2^n}+1$.
\newblock {\em Comptes Rendus Acad. Sci. Paris}, 85:329--333, 1877.

\bibitem{Schoenhage1971}
Arnold Sch{\"o}nhage and Volker Strassen.
\newblock Schnelle {Multiplikation} gro\ss er {Zahlen}.
\newblock {\em Computing}, 7:281--292, 1971.

\bibitem{Schoof1985}
Ren\'{e} Schoof.
\newblock Elliptic curves over finite fields and the computation of square
  roots $\operatorname{mod} p$.
\newblock {\em Math. Comp.}, 44(170):483--494, Apr 1985.

\bibitem{stw_phd}
Tsz-Wo Sze.
\newblock {\em On Solving Univariate Polynomial Equations over Finite Fields
  and Some Related Problems}.
\newblock PhD thesis, University of Maryland, College Park, 2007.

\bibitem{stw2011sqrt}
Tsz-Wo Sze.
\newblock On taking square roots without quadratic nonresidues over finite
  fields.
\newblock {\em Math. Comp.}, 80(275):1797--1811, July 2011.

\bibitem{lcw2008}
Lawrence~C. Washington.
\newblock {\em Elliptic Curves: Number Theory and Cryptography}.
\newblock Chapman \& Hall/CRC, 2nd edition, 2008.

\bibitem{Williams1998}
Hugh~C. Williams.
\newblock {\em {\'{E}douard Lucas and Primality Testing}}, volume~22 of {\em
  Canadian Mathematical Society Series of Monographs and Advanced Texts}.
\newblock Wiley-Interscience, 1998.

\end{thebibliography}
